\documentclass[a4paper,journal]{IEEEtran}
\usepackage[english]{babel}
\usepackage[T1]{fontenc}
\usepackage[utf8]{inputenc}
\usepackage{algorithmic}
\usepackage[square,numbers,sort&compress]{natbib}
\usepackage{graphicx}
\usepackage{pdflscape}
\usepackage[space]{grffile}
\usepackage[cmex10]{amsmath}
\usepackage{amssymb}
\usepackage{array}
\newcolumntype{H}{>{\setbox0=\hbox\bgroup}c<{\egroup}@{}}
\usepackage[caption=false,font=normalsize,labelfont=sf,textfont=sf]{subfig}
\usepackage{textcomp}
\usepackage{stfloats}
\usepackage{url}
\usepackage[allcolors=blue,bookmarks,colorlinks]{hyperref}
\usepackage{verbatim}
\usepackage{afterpage}
\usepackage{amsthm}
\theoremstyle{definition}

\usepackage{booktabs}
\usepackage[inline]{enumitem}
\usepackage{tabularx}
\newcommand*{\bbone}{\text{\usefont{U}{bbold}{m}{n}1}}
\newcommand{\Ind}{\ensuremath{\bbone}}
\providecommand{\orcidID}[1]{}

\newcommand{\zmin}{\ensuremath{z_\text{min}}}
\newcommand{\zmax}{\ensuremath{z_\text{max}}}
\newcommand{\tmax}{\ensuremath{t_\text{max}}}
\newcommand*\diff{\mathop{}\!\mathrm{d}}

\usepackage{orcidlink}

\newcommand{\customfootnote}{\footnotesize Accepted Manuscript accepted by \emph{IEEE Transactions on Evolutionary Computation} on 27 August, 2024.\newline Please cite the published version with doi:~\href{http://doi.org/10.1109/TEVC.2024.3462758}{10.1109/TEVC.2024.3462758} \copyright2024 IEEE\hfil}

\usepackage{fancyhdr}
\fancypagestyle{GlobalFootnote}{%
    \lfoot{\customfootnote} 
}
\begin{document}

\title{Using the Empirical Attainment Function for Analyzing Single-objective Black-box Optimization Algorithms}
\author{%
     Manuel López-Ibáñez\orcidlink{0000-0001-9974-1295},~\IEEEmembership{Senior Member,~IEEE},%
    \thanks{M. López-Ibáñez is with Alliance~Manchester~Business~School,~University~of~Manchester, Manchester,~M13~9PL,~UK (Email: {\tt manuel.lopez-ibanez@manchester.ac.uk}).} %
    Diederick Vermetten\orcidlink{0000-0003-3040-7162}, %
     \thanks{D. Vermetten is with Leiden Institute for Advanced Computer Science, Leiden, The Netherlands (Email: {\tt d.l.vermetten@liacs.leidenuniv.nl}).}%
  Johann Dreo\orcidlink{0000-0003-2727-9569}, %
   \thanks{J. Dreo is with Pasteur Institute, Université Paris Cité, Bioinformatics and Biostatistics Hub, Department of Computational Biology, Systems Biology Group, Paris, France (Email: {\tt johann.dreo@pasteur.fr}).}%
Carola Doerr\orcidlink{0000-0002-4981-3227}%
\thanks{C. Doerr is with Sorbonne Université, CNRS, LIP6, Paris, France (Email: {\tt carola.doerr@lip6.fr}).}%
\thanks{This is the Accepted Manuscript accepted by \emph{IEEE Transactions on Evolutionary Computation} on 27 August, 2024. Please cite the published version with doi:~\href{http://doi.org/10.1109/TEVC.2024.3462758}{10.1109/TEVC.2024.3462758}}%
}

\maketitle

\begin{abstract}
  A widely accepted way to assess the performance of iterative black-box optimizers is to analyze their empirical cumulative distribution function (ECDF) of pre-defined quality targets achieved not later than a given runtime. In this work, we consider an alternative approach, based on the empirical attainment function (EAF) and we show that the target-based ECDF is an approximation of the EAF. We argue that the EAF has several advantages over the target-based ECDF. In particular, it does not require defining a priori quality targets per function, captures performance differences more precisely, and enables the use of additional summary statistics that enrich the analysis. We also show that the average area over the convergence curves is a simpler-to-calculate, but equivalent, measure of anytime performance. 
  To facilitate the accessibility of the EAF, we integrate a module to compute it into the IOHanalyzer platform. Finally, we illustrate the use of the EAF via synthetic examples and via the data available for the BBOB suite.  
\end{abstract}


\section{Introduction}

\IEEEPARstart{I}{n} the context of benchmarking sampling-based single-objective optimization algorithms, one of the most common ways to illustrate the performance on a single problem instance is to plot the objective value of a best-so-far solution as a function of the number of solutions that have been evaluated. There are different options for the aggregation of such \emph{convergence curves} over different problem instances. In evolutionary computation, one of the most popular ways to aggregate performance data over several problems is based on the so-called \emph{empirical cumulative distribution function} (ECDF) of the number of targets attained at a given runtime. This alternative  was notably promoted by the COCO platform~\cite{HanAugMer2020coco}. To define an ECDF curve, one fixes a set of target values and counts, as a function of the solutions evaluated, how many of these targets have been achieved by the algorithm. These numbers are then normalized across all runs of the algorithm on each problem instance and then across all problem instances.   
The area under the target-based ECDF curve (AUC) can be seen as a measure of the anytime performance of the algorithm~\cite{HanAugBroTus2022anytime}.\footnote{%
  The COCO platform~\cite{HanAugMer2020coco} minimizes area over the ECDF curve, which is equivalent to the maximum runtime minus the AUC, because the former approximates the geometric average runtime.}

This \emph{target-based ECDF} generalizes the concept of time-to-target~\cite{FeoResSmi1994} or runtime distribution of a Las Vegas algorithm for a decision problem~\cite{HooStu1998uai} to multiple targets. %

One may argue that the choice of the target values for the definition of the target-based ECDF is somewhat arbitrary. We therefore consider in this paper an alternative that directly works with the probabilities that an algorithm achieves a certain quality within a given budget of function evaluations. We approximate this probability by the observed performance data, in the form of an empirical attainment function (EAF)~\cite{Grunert01}. The EAF was originally proposed for the analysis and visualization of multi-objective stochastic optimizers~\cite{GruFon2009:emaa,LopPaqStu09emaa}. In this multi-objective optimization context, the EAF estimates the probability that a single run of an algorithm attains a certain vector of objective values within a fixed runtime. %

\textbf{Our contribution:} In this paper we investigate the use of EAFs for the analysis of single-objective optimizers, previously considered in~\cite[Sec.~3.2.4]{Dreo2003phd}, \cite[Chap.~3]{ChiarandiniPhD}, \cite{Dre2009gecco}.
We first observe that the target-based ECDF is the average
value of the EAF at the pre-defined targets. We also show that the target-based ECDF with an increasing number of well-spread targets converges in the limit to an EAF-based ECDF computed as the area under the EAF across the quality dimension divided by a scaling constant. As a result, we claim that the EAF captures more information than the target-based ECDF and therefore conclude that the EAF-based ECDF has all the benefits of the target-based
approach~\cite{HanAugBroTus2022anytime} without the disadvantage of requiring an a priori set of targets. 

We also show that the AUC of the EAF-based ECDF is equivalent to the AUC of the EAF. In fact, they both measure the expected area over the convergence curve (when minimizing quality) of a single run.
Thus, any of these AUC values can be seen as a measure of anytime performance without requiring pre-defined targets.

The main drawback of the EAF is its higher computation and storage costs compared to the target-based ECDF. However, there are efficient algorithms for computing the EAF in this context and even more efficient approximations are possible. Therefore, we argue that employing
the EAF in benchmarking single-objective black-box optimizers brings
significant benefits with only a moderate overhead over the target-based ECDF. 

To facilitate the use of EAFs for the analysis of single-objective black-box optimization algorithms, we provide with this work an integration of the \texttt{eaf} package~\cite{LopPaqStu09emaa} into the IOHanalyzer platform~\cite{IOHanalyzer}. We illustrate the use of this module by analyzing performance data from the well-known Black-Box Optimization Benchmark (BBOB) collection~\cite{HanFinRosAug2009bbob}. In particular, we investigate the differences between the target-based ECDF and the EAF-based ECDF in dependence of the number of targets, and analyze how the relative ranking between algorithms may change when replacing the ECDF with the EAF as the ranking criterion.   

\textbf{Outline of the paper:} The relationship between the target-based ECDF and the EAF as well as the equivalence of various AUC metrics is discussed in Section~\ref{sec:methods}. In Section~\ref{sec:computation}, we argue that the EAF and the EAF-based ECDF can be efficiently computed in this single-objective context and we discuss their integration into IOHanalyzer. In Section~\ref{sec:IOH}, we present the results from the analysis of the BBOB data. Section~\ref{sec:stats} discusses additional statistical measures derived from the EAF and its application to the analysis  single-objective optimizers. We conclude in Section~\ref{sec:conclusions} with a summary of our arguments and perspectives for the future.

\section{Relationship between the ECDF and the EAF}\label{sec:methods}

We recall the key definitions of the target-based ECDF in Section~\ref{sec:ECDF} and those of the EAF in Section~\ref{sec:EAF}. In Section~\ref{sec:compECDFEAF} we discuss the relationship between the two notions, and in Section~\ref{sec:AUC} we treat the area under the target-based ECDF and that of the EAF-based ECDF.     

\newcommand{\dimension}{\ensuremath{D}}

\subsection{Target-based Empirical Cumulative Distribution Function}
\label{sec:ECDF}

In the following, we consider the \textit{minimization} of objective functions $f\colon \mathbb{R}^\dimension \to \mathbb{R}$, $\dimension \in \mathbb{N}^{+}$, by means of an optimization algorithm~$A$. Each run of algorithm~$A$ on a problem instance~$f$ returns a trajectory $\{(x^{(j)}, f(x^{(j)}))\mid j=1,\dotsc,\tmax\}$, with $x^{(j)}$ denoting the $j$-th point evaluated by $A$, ties broken arbitrarily in case of parallel evaluations, and $\tmax$ denoting the maximal number of function evaluations (the ``\emph{budget}''). The algorithm $A$ may be stochastic, such that different runs on the same problem instance may return different trajectories.

We focus in this work on the quality of the \textit{best-so-far solution}. To this end, we denote by 
$V(A,f,t,i):=\min\{f(x^{(j)}) \mid 1 \leq j \leq t\}$, 
the function value of the best among the first $t$ evaluated solution candidates in run $i$. 

Assuming that we have run algorithm $A$ on $f$ for $r$ independent runs, the target-based ECDF for a finite set $Z \subseteq [\zmin, \zmax] \subset \mathbb{R}$ of target values at a given budget $t$ is defined as the fraction of (run, target) pairs $(i,z)$ for which the algorithm has found, within time $t$, a solution at least as good as $z$ in the $i$-th run~\cite{HanAugBroTus2022anytime,HanAugRosFin2010comparing}.
Formally, 
\begin{equation}\label{eq:ecdf}
  \widehat{F}_Z(t) = \frac{1}{r}\sum_{i=1}^{r}\frac{1}{|Z|}\sum_{z \in Z}\Ind\left(V(A,f, t, i) \leq z\right)\enspace,
\end{equation}
where  $\Ind({\mathcal{C}})$ denotes the indicator function that is $1$ when the condition $\mathcal{C}$ is satisfied and $0$ otherwise. 

Other definitions of ECDF may be found in the literature. For example, the target-based ECDF defined above is not the same as the classical multivariate ECDF~\cite{GruFon2002spl}.

\subsection{Empirical Attainment Function}
\label{sec:EAF}

In the context of benchmarking the anytime performance of single-objective optimizers, given a pair  $(t,z)$ of budget and function value, respectively, we say that run $i$ of algorithm $A$ \emph{attains} $(t, z)$ if and only if algorithm $A$ solving problem $f$ obtains an objective function value not worse than $z$ not later than $t$, that is, if and only if $V(A,f, t, i) \leq z$. We can formulate the EAF in this context as the fraction of $r$ runs of algorithm $A$ on function $f$ that attain an objective value $z$ at budget $t$:
\begin{equation}\label{eq:eaf}
\widehat{\alpha}(t,z) = \frac{1}{r}\sum_{i=1}^{r}\Ind\left(V(A,f, t, i) \leq z\right)\enspace.
\end{equation}

The above matches the original formulation of the EAF~\cite{Grunert01}, if we consider that each run $i$ produces a set $S_i  \subset \mathbb{R}^2$ of points 
$\left(t, V(A,f,t,i)\right)$, i.e., \mbox{$S_i = \{ \left(t, V(A,f,t,i)\right) \mid 1 \leq t \leq \tmax\}$}.
In this original formulation, a set \mbox{$S_i \subset \mathbb{R}^d$} \emph{attains} the point $y=(y_1,\dotsc, y_d) \in \mathbb{R}^d$ iff there exists at least one point $s=(s_1,\dotsc,s_d) \in S_i$ such that, for all $1 \le j \le d$, it holds that $s_j \le y_j$. In this case, the value of the EAF of a set of sets $\{S_1, \dotsc, S_n\}$ at a given point $y \in  \mathbb{R}^d$ is defined as the fraction of sets $S_i$, $i=1, \dotsc, n$, that attain that point $y$.
See Figure~\ref{fig:EAF_simple} for an illustrative example of the EAF for three sets, corresponding to three runs of an algorithm.

\subsection{Comparison of the target-based ECDF and the EAF}
\label{sec:compECDFEAF}

If we have a finite number of targets $|Z| \in \mathbb{N}^{{+}}$, the target-based ECDF value for a given budget $t$ \eqref{eq:ecdf} equals the average EAF value in these targets for time $t$, as the following simple swap of summands shows. 
\begin{equation}
\begin{split}
\frac{1}{|Z|}\sum_{z \in Z}\widehat{\alpha}(t,z) \\
= \frac{1}{|Z|}\sum_{z \in Z}\left(\frac{1}{r}\sum_{i=1}^{r}\Ind \left(V(A,f, t, i) \leq z\right)\right)\\
= \frac{1}{r}\sum_{i=1}^{r}\frac{1}{|Z|}\sum_{z \in Z}\Ind \left(V(A,f, t, i) \leq z\right) = \widehat{F}_Z(t)\enspace.
\end{split}
\end{equation}

\newcommand{\eafECDF}{\ensuremath{F_{\widehat{\alpha}}(t)}}

With the above equality at hand, it is easy to show that the target-based ECDF of a given $t$, evaluated for increasingly finer, well-spread partitions of $[\zmin,\zmax]$ converges to the following partial integral
\begin{equation}\label{eq:eaf_ecdf}
  \eafECDF:=\frac{1}{\zmax- \zmin} \int_{\zmin}^{\zmax} \widehat{\alpha}(t,z)\diff z
\end{equation}
of the EAF over the quality dimension $z$ in this interval, evaluated for the same budget $t$. That is, the target-based ECDF with infinite number of well-spread targets is the area under the EAF at a fixed value of $t$ divided by a scaling constant. Thus, the integral $\eafECDF$ 
  defines an \emph{EAF-based ECDF} that serves the same purpose as the target-based ECDF and does not require predefined targets.

Formally, let \mbox{$Z=\{z_1,\dots,z_{|Z|} \}\subset[\zmin,\zmax]$} with \mbox{$z_i < z_{i+1}$} for all \mbox{$1 \le i \le |Z|-1 $}, and let 
$\delta(Z)$ $:=$ $\max\left\{z_1-\zmin,\right.$ $\max_{1 \le i \le |Z| -1} \{z_{i+1} -z_i\},$ $\left.\zmax-z_{|Z|}  \right\}$ 
be the largest ``gap'' between any two neighboring target values in the set $Z \cup\{\zmin,\zmax\}$.  
Then, one can show that
\begin{equation}
\eafECDF = \lim_{\delta(Z) \to 0} \widehat{F}_Z(t)
\end{equation}
in the same way as the Riemann sum of a function converges to the integral of the function:
\begin{equation}\label{eq:approx}
\begin{split}
\lim_{\delta(Z) \to 0} \widehat{F}_Z(t) \\
= \frac{1}{r}\sum_{i=1}^{r}\frac{1}{\zmax- \zmin}\int_{\zmin}^{\zmax}\Ind\left( V(A,f, t, i) \leq z \right) \diff z\\
= \frac{1}{\zmax- \zmin} \int_{\zmin}^{\zmax}\frac{1}{r}\sum_{i=1}^{r}\Ind \left(V(A,f, t, i) \leq z \right) \diff z\\
=\frac{1}{\zmax- \zmin} \int_{\zmin}^{\zmax} \widehat{\alpha}(t,z)\diff z \enspace = \eafECDF\enspace.
\end{split}
\end{equation}

\begin{figure}
  \begin{minipage}{1.0\linewidth}
    \centering%
    \includegraphics[width=.8\linewidth]{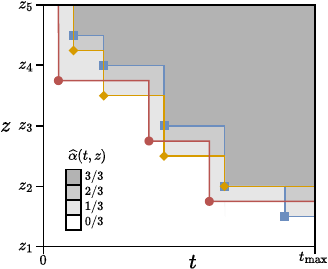}
    \caption{Visualization of the EAF corresponding to three runs (red circles, orange diamonds, and blue squares) of an hypothetical single-objective optimizer. For each run, we record as a point, the  runtime $t$ (x-axis) at which a new best-so-far objective function value $z$ (y-axis) was found (assuming minimization). The value of the EAF $\widehat{\alpha}$ at point $(t, z)$ is given by the shade of gray specified in the legend. We mark 5 quality targets in the y-axis $\{z_1, z_2, \dots, z_5\}$ for computing the target-based ECDF in Figure~\ref{fig:EAF_ECDF}.}
    \label{fig:EAF_simple}
  \end{minipage}
  \vskip 1.5em
  \begin{minipage}{1.0\linewidth}
    \centering%
    \includegraphics[width=.8\linewidth]{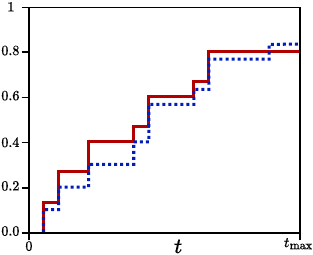}
    \caption{Target-based ECDF (solid red line) and EAF-based ECDF (blue dashed line) corresponding to the runs shown in Figure~\ref{fig:EAF_simple}. The target-based ECDF uses the quality targets $\{z_1, z_2, \dots, z_5\}$ defined in Figure~\ref{fig:EAF_simple}. For the EAF-based ECDF, we only need $\zmin = z_1$ and $\zmax = z_5$. 
}
     \label{fig:EAF_ECDF}
  \end{minipage}
\end{figure}

Figure~\ref{fig:EAF_ECDF} shows the target-based ECDF and the EAF-based ECDF corresponding to the example runs shown in Figure~\ref{fig:EAF_simple}. This example illustrates that the target-based ECDF with a small number of targets may either over-estimate or under-estimate the EAF-based ECDF, an observation that will be further discussed in Section~\ref{sec:sensitivity}.

\subsection{Area under the EAF as a measure of anytime performance}\label{sec:AUC}

The area under the curve (AUC) of the target-based ECDF has been used in the literature as a measure of the anytime performance of single-objective optimizers~\cite{HanAugBroTus2022anytime,YeDoeWang2022ac}. As shown above, the target-based ECDF is an approximation to the EAF-based ECDF~\eqref{eq:eaf_ecdf}. Thus, the AUC of the EAF-based ECDF would be a more precise measure of anytime performance. The example in Fig.~\ref{fig:EAF_ECDF} already shows that the AUC values of the two curves are not equivalent. 

The AUC of the EAF-based ECDF is given by:
\begin{equation}\label{eq:auc_eaf_ecdf}
  \text{AUC}\left(\eafECDF\right) = \frac{1}{\zmax-\zmin} \int_{1}^{\tmax} \diff t\int_{\zmin}^{\zmax} \widehat{\alpha}(t,z)\diff z\enspace.
\end{equation}

By definition, the double integral in the AUC of the EAF-based ECDF \eqref{eq:auc_eaf_ecdf} is equivalent to  the area under the EAF given by the following Lebesgue integral:
\begin{equation}\label{eq:auc_eaf}
\text{AUC}\left(\widehat{\alpha}(t,z)\right) = \int_{(1,-\infty)}^{(\tmax, \zmax)} \widehat{\alpha}(t,z)\diff t\diff z\enspace,
\end{equation}
which is equivalent to computing the well-known hypervolume metric in
multi-objective optimization~\cite{ZitThi1998ppsn} of the three-dimensional points
$\{(t,z, 1-\widehat{\alpha}(t,z))\mid t \in [1, \tmax], z \in [-\infty,\zmax]\}$ with reference point $(\tmax, \zmax, 1)$. This computation only requires $\Theta(n \log n)$~\cite{BeuFonLopPaqVah09:tec}, where $n$ is the number of points.

The only differences between \eqref{eq:auc_eaf_ecdf} and \eqref{eq:auc_eaf} are the normalization constant $1/(\zmax-\zmin)$ and the presence of finite $\zmin$ in the integral. If we wish to define $\zmin$ in \eqref{eq:auc_eaf}, e.g., for normalization purposes, we only have to clip $V(A,f,t,i)$ so that its minimum is never lower than $\zmin$. Evidently, if $V(A,f,t,i)$ is never lower than $\zmin$, e.g., because $\zmin$ is the optimal value of $f$ or better, then $\zmin$ has no effect in the AUC values, except for the normalization constant. With this in mind, the presence of a finite $\zmin$ is not necessary for calculating the AUC of the EAF \eqref{eq:auc_eaf}, thus we use $-\infty$ for generality.%

Another measure of anytime performance proposed in the literature~\cite{LopStu2013ejor} is the area over the convergence curve (AOCC),
$S_i = \{ (t, V(A,f,t,i)) \mid 1\leq t\leq \tmax \}$, which can be calculated by summing up the areas of the rectangles formed by the points $(t, V(A,f,t,i))$ and $(t+1, \zmax)$ as follows:
\begin{equation}\label{eq:aoc}
\text{AOCC}(S_i) = \sum_{t=1}^{\tmax - 1} \max\{0, \zmax - V(A,f,t,i)\}\enspace.
\end{equation}

We can see that Eq.~\eqref{eq:auc_eaf} is equivalent to the mean AOCC value of the multiple runs if we expand the former as follows:
\begin{equation}\label{eq:auc_eaf_proof_1}
\begin{split}
     \int_{1}^{\tmax} \diff t\int_{-\infty}^{\zmax} \widehat{\alpha}(t,z)\diff z\\
    =    \int_{1}^{\tmax} \diff t\int_{-\infty}^{\zmax} \frac{1}{r} \sum_{i=1}^{r}\Ind\left(V(A,f, t, i) \leq z\right)\diff z\\
    =  \frac{1}{r} \sum_{i=1}^{r} \int_{1}^{\tmax} \diff t\int_{-\infty}^{\zmax}  \Ind\left(V(A,f, t, i) \leq z\right)\diff z\enspace.
  \end{split}
\end{equation}

For a fixed $V(A,f, t, i)$, the indicator function $\Ind\left(V(A,f, t, i) \leq z\right)$ is $1$ from \mbox{$z=V(A,f, t, i)$ up to $\infty$}, and $0$ otherwise.   Thus, its integral over $z$ bounded by $(-\infty, \zmax]$ is $\zmax - V(A,f,t,i)$ when $\zmax \geq V(A,f,t,i)$ and zero otherwise. Therefore,
\begin{equation}\label{eq:auc_eaf_proof_2}
\begin{split}
  \frac{1}{r} \sum_{i=1}^{r} \int_{1}^{\tmax} \diff t\int_{-\infty}^{\zmax}  \Ind\left(V(A,f, t, i) \leq z\right)\diff z\\
  = \frac{1}{r}\sum_{i=1}^r \sum_{t=1}^{\tmax - 1} \max\{0, \zmax - V(A,f,t,i)\}\\
    = \frac{1}{r}\sum_{i=1}^r \text{AOCC}(S_i)\enspace.
  \end{split}
\end{equation}

In other words, not only the AUC of the EAF-based ECDF is equivalent to the
AUC of the EAF (multiplied by a normalization constant), but also both are equivalent to the mean AOCC of the runs used to compute the EAF. A more general proof of the relationship between the expected value of a generalized hypervolume indicator (i.e., where the AOCC is multiplied by a weight function) and the first-order attainment function (i.e., the EAF when the number of runs goes to infinity) is provided by \cite{GruFon2012ea}. The simpler proof in Eq.~\eqref{eq:auc_eaf_proof_1}--\eqref{eq:auc_eaf_proof_2} does not consider the weighted variant of the hypervolume indicator, but it applies to a finite number of empirical samples.

There are situations in which we are only interested in the AUC value, for example, for benchmarking purposes or when using the AUC as the metric that guides the tuning of parameters. In those situations, instead of storing the data for all runs to compute the target-based ECDF or the EAF and their AUC value, it is much more convenient to compute the AOCC of each run as soon as it finishes and only store the AOCC values. Moreover, recomputing the mean AOCC after performing additional runs is trivial, whereas updating the EAF or the ECDF is less so. The use of the mean AOCC to guide the tuning of parameters was already demonstrated in the literature~\cite{LopStu2013ejor}, but the connection with the AUC of the target-based ECDF was not known at the time.

\section{Computation of the EAF}
\label{sec:computation}
\subsection{Efficiency of computation}

In the formulation of the EAF discussed here for benchmarking single-objective black-box optimization algorithms~\eqref{eq:eaf}, points are pairs $(t,z)$ %
of runtime and objective function value. The EAF of an algorithm at this point is the probability that the algorithm has found, within the first $t$ function evaluations, a solution that is at least as good as~$z$. 
In this two-dimensional context, the computation of the EAF requires $\Theta(n \log n + rn)$ time/steps~\cite{FonGueLopPaq2011emo}, where $r$ is the number of runs and $n \le r \tmax$ is the total number of points, i.e., objective value improvements recorded across all runs. 

The number of points required to fully represent the EAF is  $\Theta(rn)$~\cite{FonGueLopPaq2011emo}. Reducing $n$ by merging points representing very quick or very small improvements in succession is an easy way to speed-up the computation of the EAF and reduce its size. 
Therefore, we argue that computing or storing the EAF is feasible for a reasonable number of runs and points.

The computation of the EAF-based ECDF~\eqref{eq:eaf_ecdf} may seem more complicated than the target-based ECDF~\eqref{eq:ecdf}. However,
given that $\widehat{\alpha}(t,z)$ for a fixed $t$ is a step function monotonically increasing with respect to $z$, the integral $\eafECDF$ can be seen as the Lebesgue integral over $z$ for each $t$, which is easy to calculate as a sum of rectangles. Moreover, we need $\zmin$ to show the exact equivalence between the target-based and EAF-based ECDFs, however, the computation of the EAF does not require defining $\zmin$.

\subsection{Integration into IOHanalyzer}

For practical usage of the EAF with various types of benchmarking data, we have integrated several ways of using EAF-based analysis into IOHanalyzer~\cite{IOHanalyzer}, which is a part of the IOHprofiler benchmarking platform. For this purpose, the computations of the EAF are calculated by version 2.4.1 of the \texttt{eaf} R package~\cite{LopPaqStu09emaa}.\footnote{\url{https://mlopez-ibanez.github.io/eaf/}}  Through this integration, several ways of computing, visualizing and comparing EAF-based statistics, are available directly in the web-interface of IOHanalyzer (and the corresponding R package). Benchmarking data from platforms such as COCO~\cite{HanAugMer2020coco}, Nevergrad~\cite{nevergrad}, and IOHexperimenter~\cite{IOHexperimenter2024} can also be analyzed using these methods.

To speed up computations when creating visualizations, we make use of a subsampling approach when calculating the EAF, ECDF, and AUC values. This subsampling is done on the runtime-axis, where we select a set of $50$ log-spaced runtimes at which we calculate the fraction of targets hit for the target-based ECDF
or at which we take the function values to pass to the EAF calculation.
Given this subsampling, the computation of the EAF and the EAF-based ECDF took just 432 seconds in a single CPU for all datasets considered in the next section (1295 datasets in total), which is a reasonable time for the analysis of an experiment of this scale.

\section{Use in Practice}
\label{sec:IOH}

We illustrate the use of the EAF for benchmarking single-objective black-box optimizers by analyzing data from two commonly used optimization algorithms: BFGS~\cite{Bro1970bfgs,Fle1970bfgs,Gol1970bfgs,Sha1970bfgs}
 (from the \texttt{scipy} package~\cite{Scipy2020natmet}) and CMA-ES~\cite{HanOst1996cma} (from the \texttt{pyCMA} package~\cite{HanAkiBau2019pycma}). %
Both these algorithms are run on the 24 functions from the BBOB suite. All functions allow specifying the dimension of the decision space $\dimension$.
Each algorithm is run 15 times on each function with different random seeds, each time using a different instance (i.e., rotation and translation) of the function. Each individual run is stopped after $\tmax=10^4\cdot \dimension$ function evaluations. %

\begin{figure}
  \centering
  \begin{minipage}{1.0\linewidth}
\includegraphics[width=\textwidth]{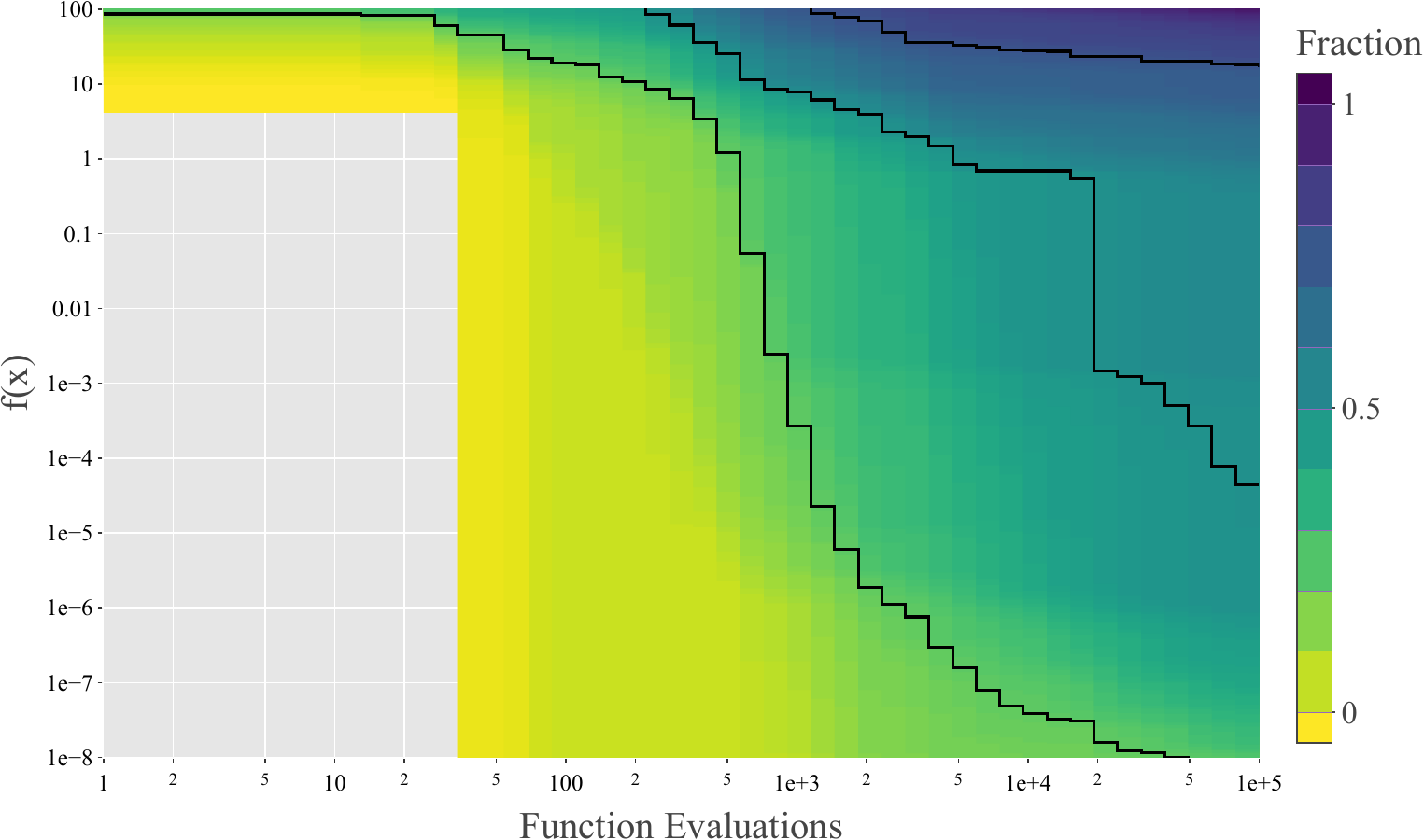}
\caption{EAF of BFGS over the 24 10-dimensional BBOB functions. The color of a point gives the fraction of runs that reach a given value of $f(x)$ not later than a given number of function evaluations. The black lines delimit the regions with a value (from bottom to top) $\leq 0.25$, $\leq 0.5$ and $\leq 0.75$.}\label{fig:eaf_bfgs}
\end{minipage}
\\[1em]
  \begin{minipage}{1.0\linewidth}
    \centering
\includegraphics[width=\textwidth]{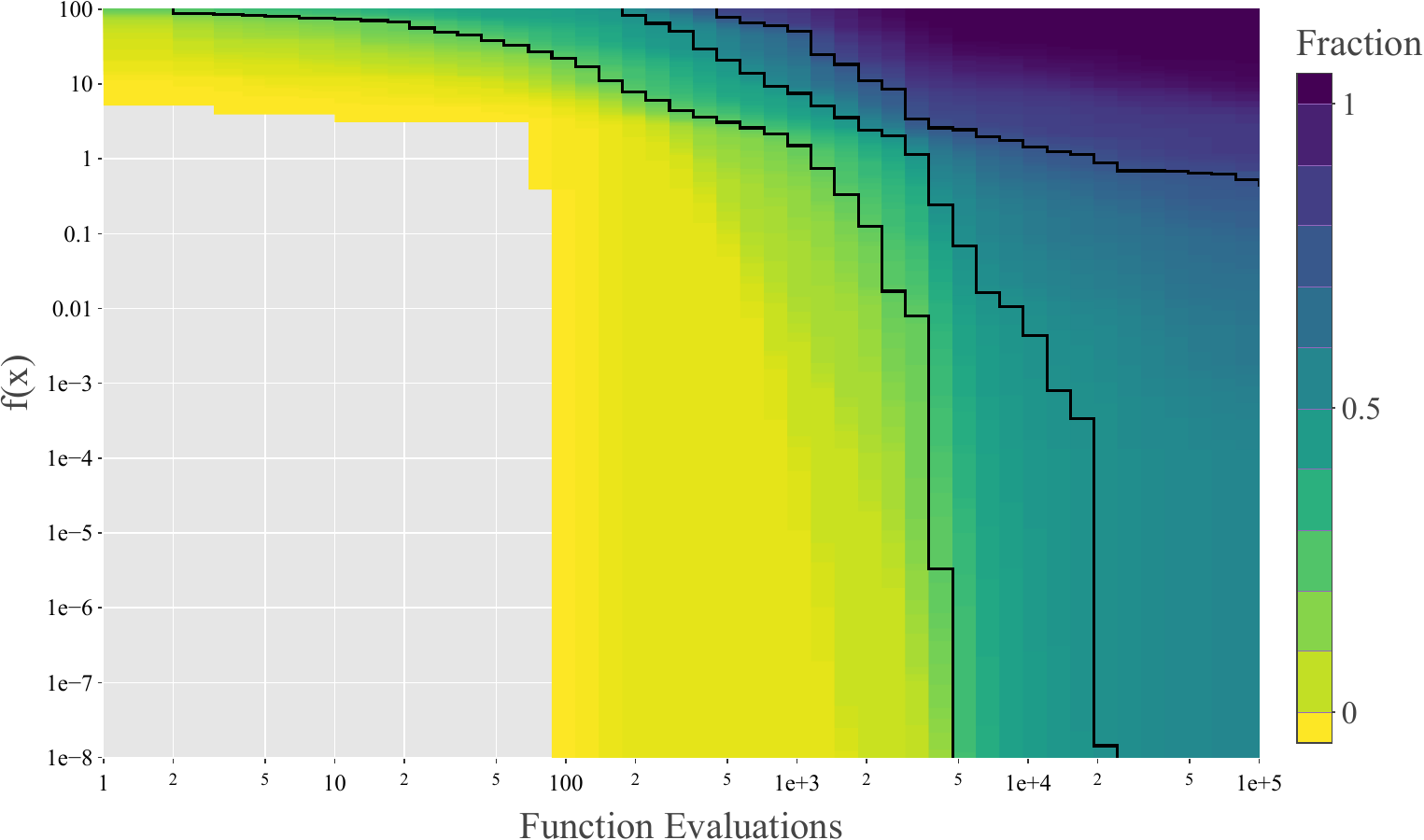}
    \caption{EAF of CMA-ES over the 24 10-dimensional BBOB functions. See the caption of Figure~\ref{fig:eaf_bfgs} for more details.}
    \label{fig:eaf_cma}
  \end{minipage}
\end{figure}

For computing the EAF-based ECDF \eqref{eq:eaf_ecdf} and its AUC on the BBOB functions, we follow the convention set by COCO~\cite{HanAugMer2020coco} and define $\zmin = 10^{-8}$ and $\zmax = 10^2$. To match the logarithmic scale of the targets used by COCO, we use a logarithmic scaling of the interval $[\zmin, \zmax]$ (effectively taking the log-precision as our function values throughout any EAF-based ECDF calculations).  Figures~\ref{fig:eaf_bfgs} and \ref{fig:eaf_cma} show the resulting EAF of the BFGS and CMA-ES algorithms, respectively, aggregated over all 24 BBOB functions in dimension 10. Darker colors indicate a higher probability of attaining a certain objective function value within the number of function evaluations given by the x-axis. From these figures, we see a clear trend in the optimization behavior, where BFGS converges faster at the beginning of the search (darker colors on the lower left of the colored region), but quickly stagnates (small contrasts going from left to right). We highlight specific levels of the EAFs $(0.25, 0.5, 0.75)$ as black lines, which suggest that the results of CMA-ES are more robust because the $0.5$ and $0.75$ lines of CMA-ES almost completely dominate the corresponding ones of BFGS.  %

\begin{figure}
    \centering
    \includegraphics[width=0.5\textwidth]{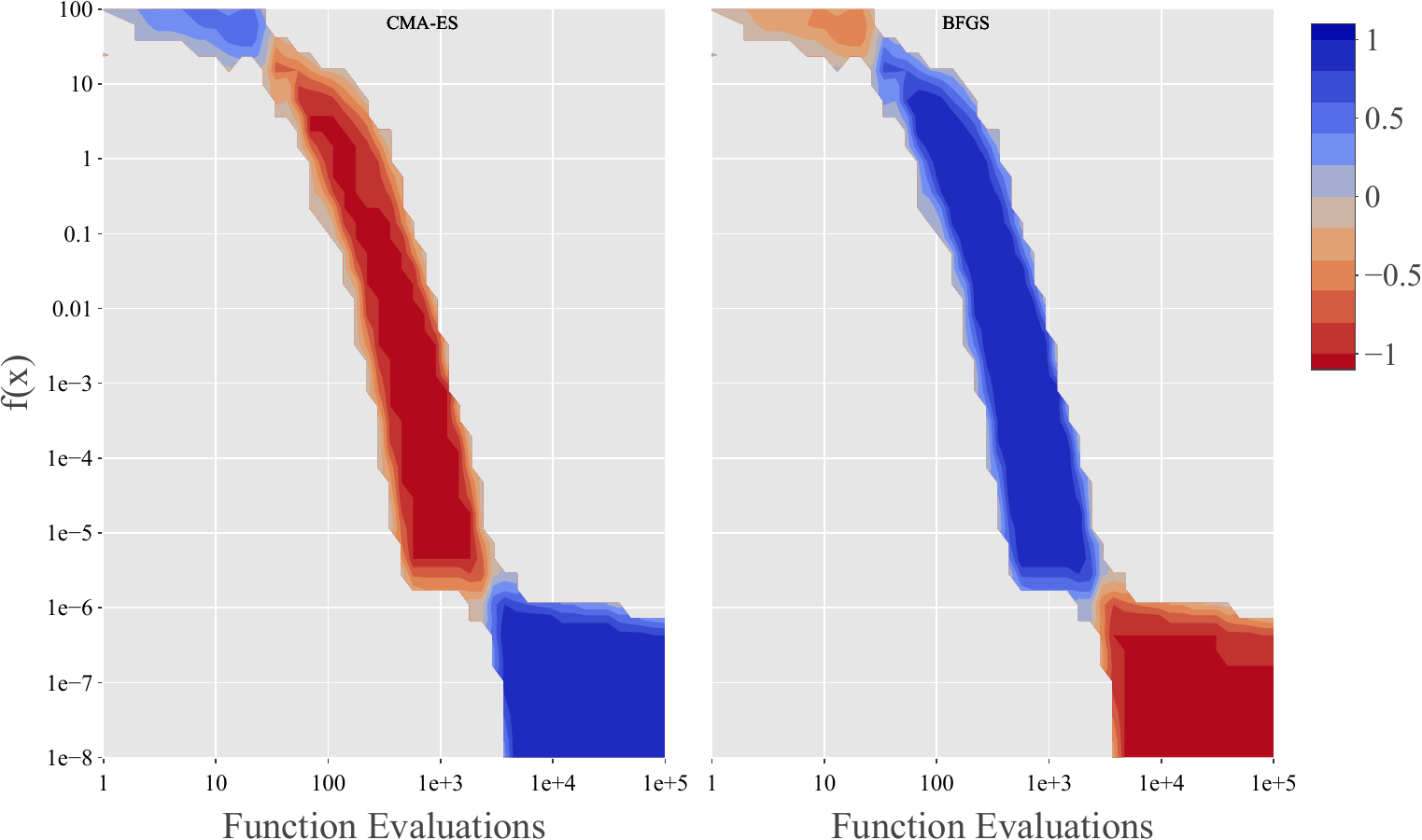}
    \caption{Differences in the EAFs of CMA-ES and BFGS 
    on BBOB function 14, in 10 dimensions. Blue colors indicate regions where the algorithm in the plot title outperforms the other algorithm.}
    \label{fig:eafdiff_bfgs_cma}
\end{figure}

\begin{figure*}
    \centering
    \includegraphics[width=.9\linewidth]{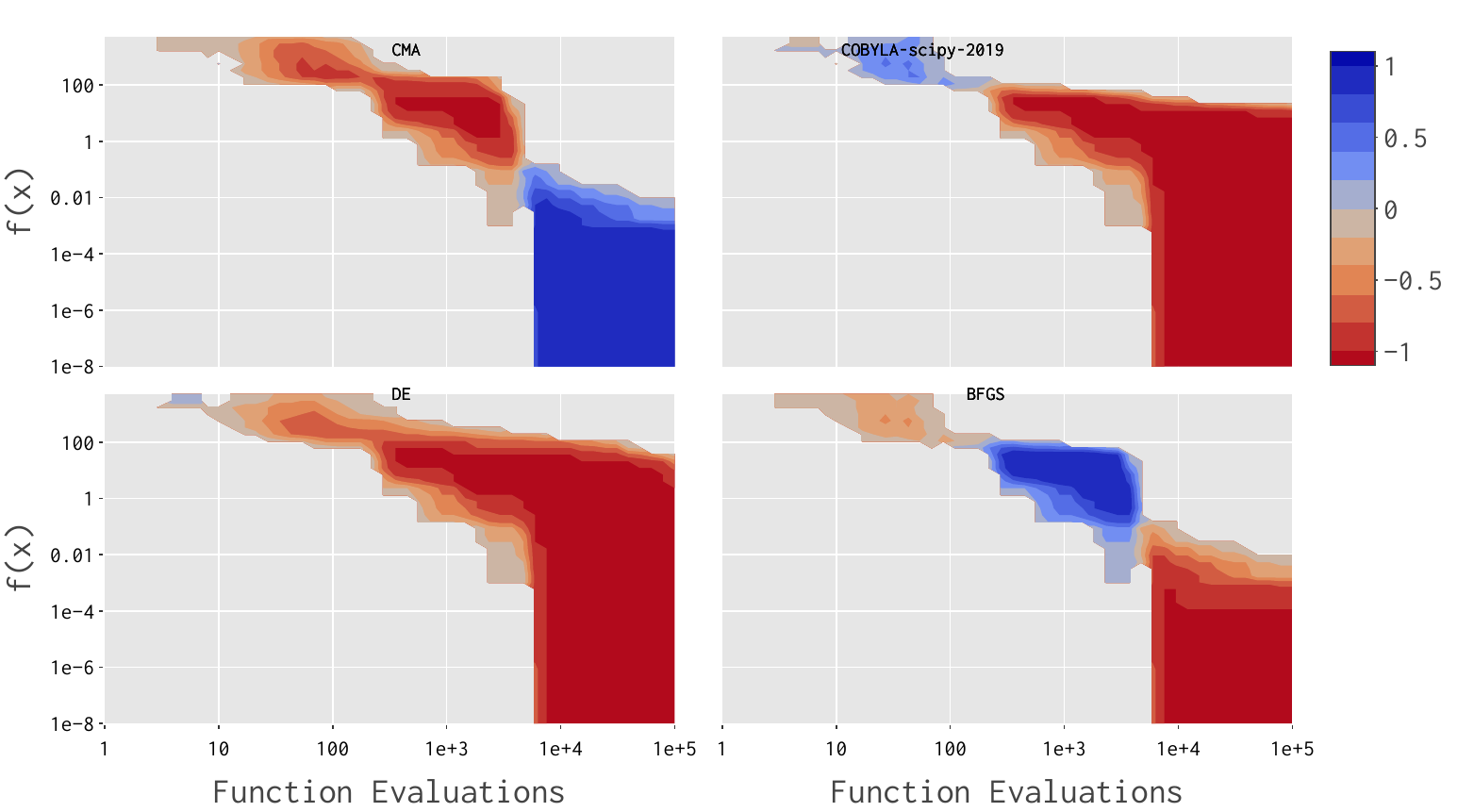}
    \caption{Differences in the EAFs of a portfolio of 4 algorithms on BBOB function 11, in 10 dimensions. Blue colors indicate regions where the algorithm in the plot title outperforms the other algorithms.}
    \label{fig:eafdiff_portfolio}
\end{figure*}

In addition to the individual EAF plots, we can also perform a comparative analysis by plotting the differences between two EAFs. This is illustrated in Figure~\ref{fig:eafdiff_bfgs_cma}, where the differences between CMA-ES and BFGS on a single BBOB function are shown. In this example, BFGS outperforms CMA-ES between $50$ and $1\text{e}{+}3.5$ (approx.~$3162$)
evaluations, whereas CMA-ES wins after that, managing to converge to $f(x)=10^{-8}$ more reliably. %
This analysis can additionally be extended to highlight the difference between an algorithm and a portfolio, by taking the upper envelope of the portfolio's EAF. This portfolio-based comparison is illustrated in Figure~\ref{fig:eafdiff_portfolio}, where a set of four algorithms are compared on benchmark function F11. The selected portfolio consists of CMA-ES, BFGS, Differential Evolution~\cite{StoPri1997:de} (DE) and Cobyla~\cite{Powell1994cobyla}, and their performance data is taken from the COCO data repository. From this figure, we see that in three stages of the search, Cobyla, BFGS, and CMA-ES, respectively, outperform every other considered algorithm, whereas there is no pair $(t,z)$ with $t>10$ in which the EAF of DE
dominates the maximum of the EAF values of the other three algorithms. %

\subsection{Sensitivity of the ECDF wrt Target Values}\label{sec:sensitivity}

\begin{figure}
    \centering
    \includegraphics[width=\linewidth]{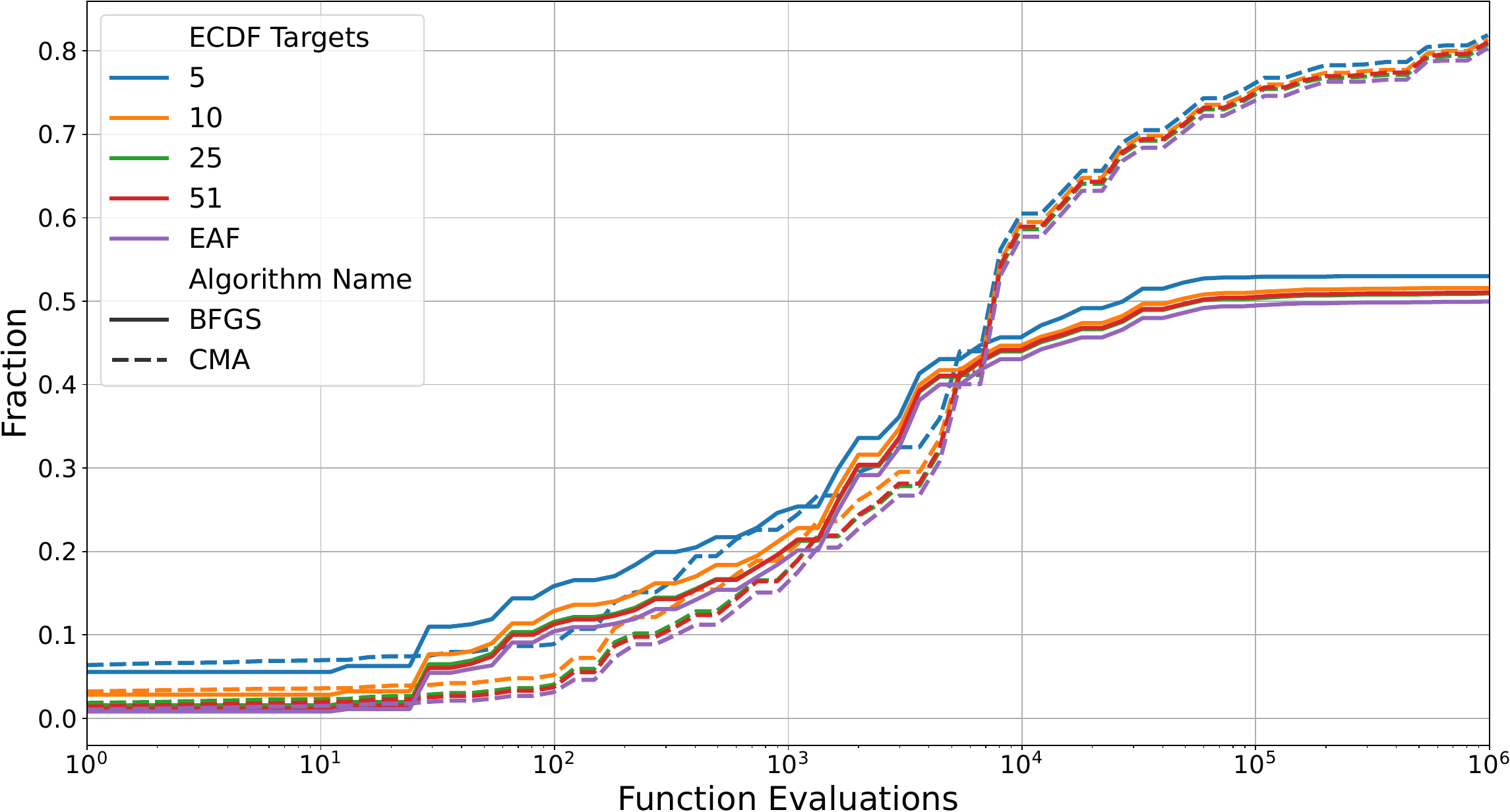}
    \caption{ECDF as calculated using the EAF \eqref{eq:eaf_ecdf} and the target-based approach \eqref{eq:ecdf}, with different numbers of targets, for two different algorithms, BFGS and CMA-ES. Data is aggregated over all 24 10-dimensional BBOB functions.}
    \label{fig:ecdf_target_impact}
\end{figure}

In this subsection, we study the sensitivity of the target-based ECDF to the number of pre-defined target values. That is, we empirically evaluate the convergence speed of the limit in~\eqref{eq:eaf_ecdf}.
Following COCO, we use a logarithmic scale by selecting equally spaced targets that are powers of 10 between $10^{-8}$ and $10^{2}$. In addition to the default of COCO with 51 targets, we also consider 5, 10, and 25 targets. For each of these values and each algorithm under comparison, we compute the target-based ECDF \eqref{eq:ecdf} and its difference in value to the EAF-based ECDF \eqref{eq:eaf_ecdf}. To increase the amount of algorithms being compared, we now make use of the data from the BBOB repository that matches the setup from the previous section: 1 run on each of 15 different instances. For the 10-dimensional functions, this leads to a set of 211 algorithms.\footnote{For readers familiar with the COCO platform, we note that we do not apply the bootstrapping approach suggested in~\cite{HanAugMer2020coco} to simulate restarts. Instead, we compute the EAF from the available trajectories.} 

In Figure~\ref{fig:ecdf_target_impact}, we illustrate the impact of changing the number of targets for the BFGS and CMA-ES. As can be seen from this figure, the ECDF based on a limited number of targets consistently overestimates the EAF-based ECDF.

\begin{figure}
    \centering
    \includegraphics[width=0.5\textwidth]{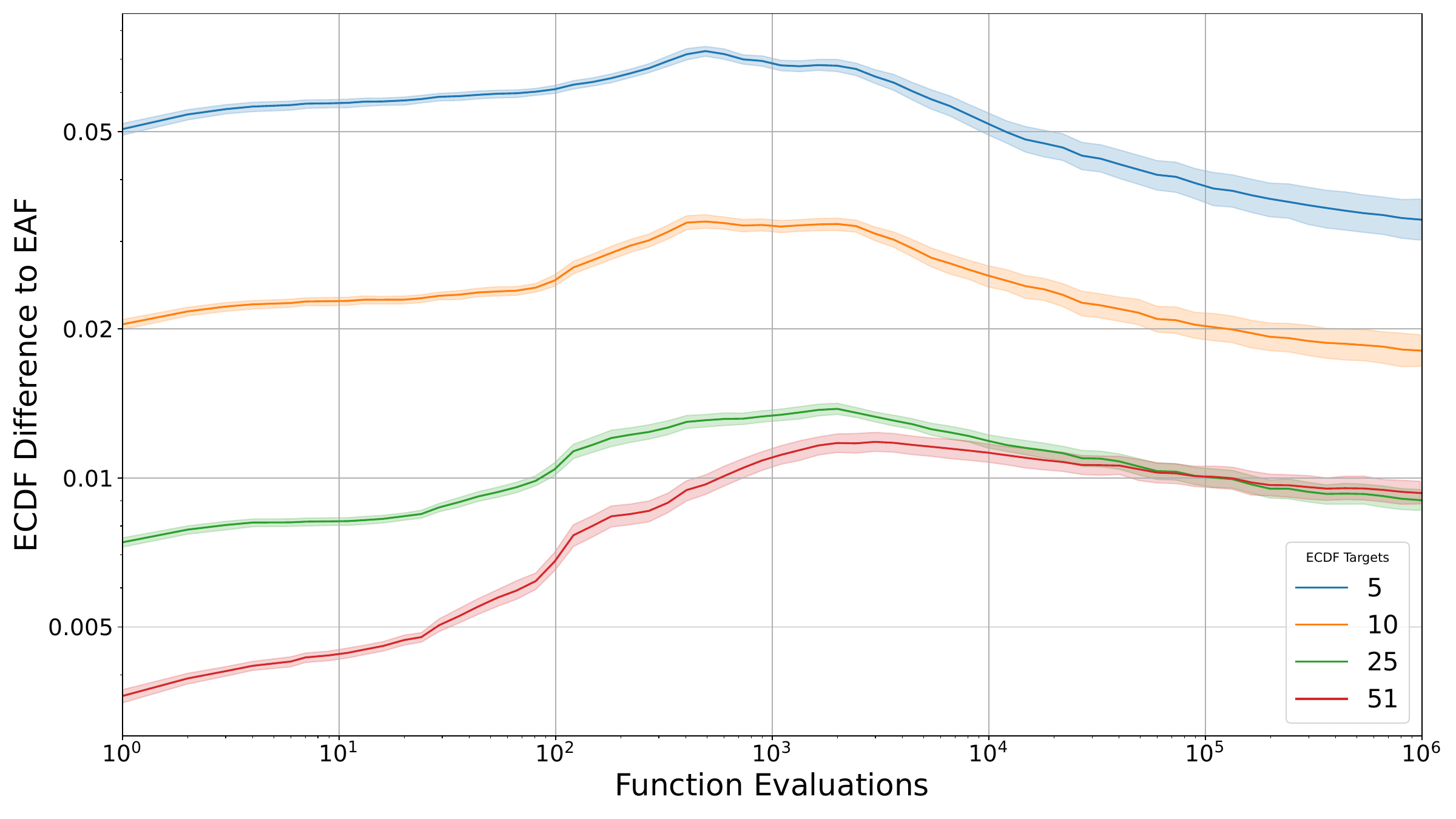}
    \caption{Mean and standard deviation of differences in target-based ECDFs with different numbers of targets relative to the EAF-based ECDF. Each individual ECDF value is an aggregation of the 24 10-dimensional BBOB problems for one algorithm. These are aggregated over all $211$ processed algorithms from COCO's repository. Both axes are in logarithmic scale.}%
    \label{fig:ecdf_diffs_targets}
\end{figure}

Figure~\ref{fig:ecdf_diffs_targets} shows the mean and standard deviation of
these differences over all 211 algorithms and for the different numbers of
targets. As expected, the differences decrease rapidly with increasing number
of targets. For the data in this experiment and with 51 targets, differences %
are larger for higher
runtime. In addition, we clearly see that the overestimation of the
target-based ECDF was not an artifact of the selected algorithms in
Figure~\ref{fig:ecdf_target_impact}, but a consistent phenomenon.  The
differences between the target-based and EAF-based ECDF values are a
combination of two distinct effects. 
The first and more important effect is inherent in the approximation induced
by the choice of targets: any objective value attained between two successive
targets will result in the same value of the target-based ECDF, whereas any
objective value improvement increases the value of the EAF-based ECDF. For a
single run, the size of the differences caused by this effect is at most
$1/|Z|$, where $|Z|$ is the number of targets.

A secondary effect is due to the difference between $\zmax$ and the first
target. If $\zmax$ equals the first target, then the ECDF value when attaining
$\zmax$ differs between the two approaches. For simplicity, let us consider a
single run that exactly attains $\zmax$ at some $t$. In the target-based ECDF,
this run has a value of $1/|Z|$, while this same run has value $0$ in the
EAF-based ECDF, at time $t$. This initial difference causes the target-based ECDF to
consistently overestimate the EAF-based ECDF, as seen in
Figure~\ref{fig:ecdf_diffs_targets}. We could remove this initial difference by
increasing $\zmax$ relative to the first target by a factor of
$(\zmax - \zmin)/(|Z|-1)$ before computing the EAF. However, this adjustment
would actually cause under-estimation due to the first effect mentioned
above. Moreover, a run that attains an objective value strictly between the
adjusted $\zmax$ and the first target will result in a positive value of the
EAF-based ECDF, but a value of $0$ of the target-based ECDF, thus also leading
to under-estimation. %
In other words, if $\zmax$ is larger than the first target, then the
target-based ECDF may initially under-estimate the EAF-based ECDF. In summary,
due to the unavoidable first effect, adjusting $\zmax$ to remove the secondary
effect does not eliminate the differences, it just shifts them from
consistent over-estimation to primarily under-estimation. %
In addition, applying this adjustment would require computing an EAF-based ECDF
for each set of targets, which greatly complicates the comparison between
target-based and EAF-based ECDFs. Thus, we have decided to define $\zmax$ as
the first target in all our experiments.

\begin{figure}
    \centering
    \includegraphics[width=\linewidth]{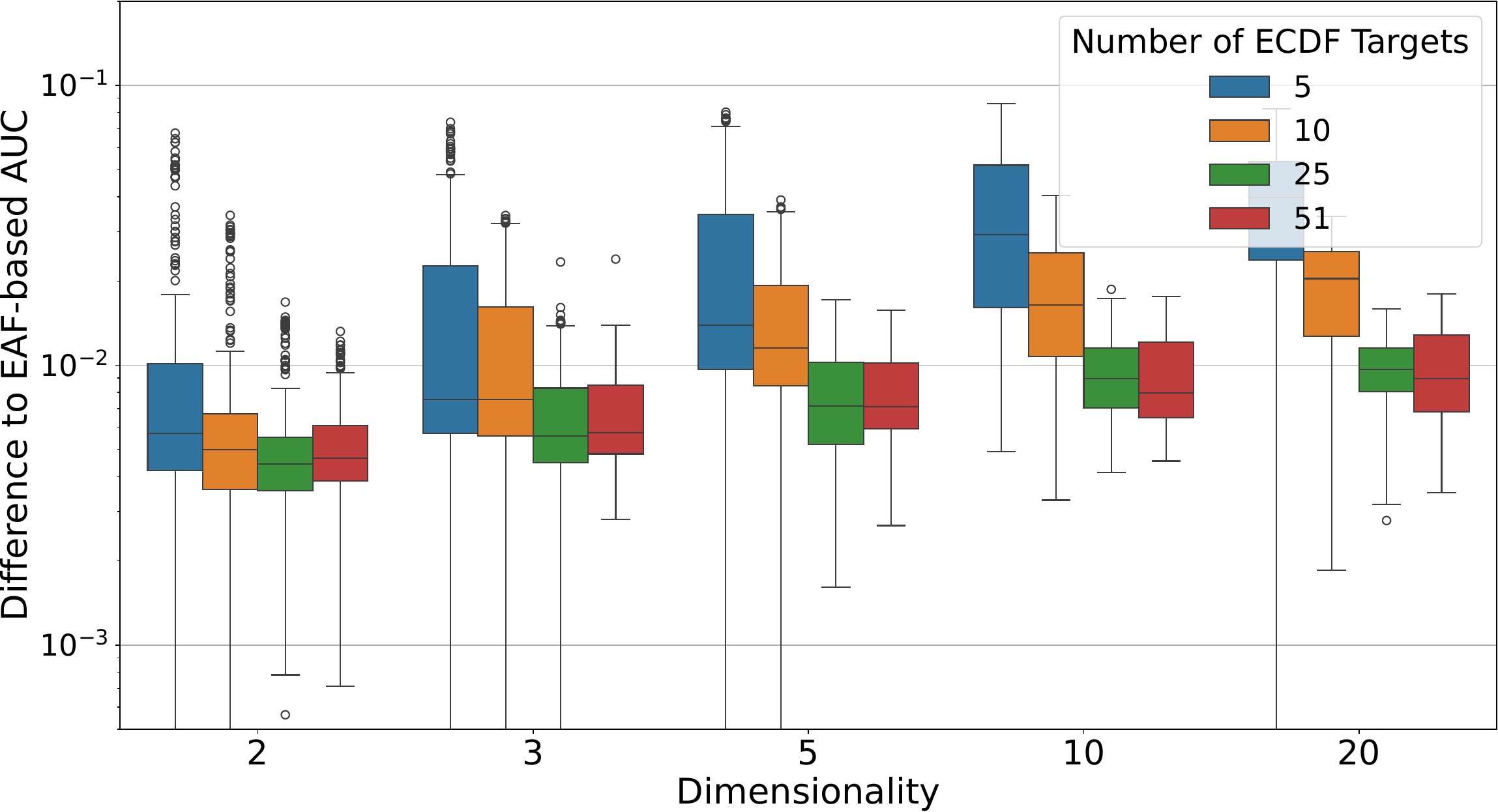}
    \caption{Distribution of differences between the AUC of the EAF-based ECDF and the AUC of the target-based ECDF, for different number of targets. Each data point is the mean AUC difference for a single algorithm averaged over 15 runs on the 24 BBOB problems of the dimension given in the x-axis %
    and the box plots summarize the data distribution of all 211 processed algorithms from the COCO repository}
    \label{fig:auc_diffs}
\end{figure}

In Figure~\ref{fig:auc_diffs}, we can see that the differences observed earlier
also translate into differences in AUC. We can again observe that the
differences between the AUCs of the target-based ECDF and EAF-based ECDF
decrease quickly with the number of targets. At the same time, the effect of
the number of targets is stronger for problems in higher dimensions. We can also
observe that going from 25 to 51 targets does not help much in reducing the
differences below $1\%$. The relatively large variances of the boxplots
indicate that the low variance observed in Figure~\ref{fig:ecdf_diffs_targets}
does not translate into consistent differences across algorithms. In other
words, the estimation error of the target-based ECDF is larger for some algorithms.

\subsection{Changes in ranking when using ECDF vs EAF}

\begin{figure*}
    \centering
    \includegraphics[width=\textwidth]{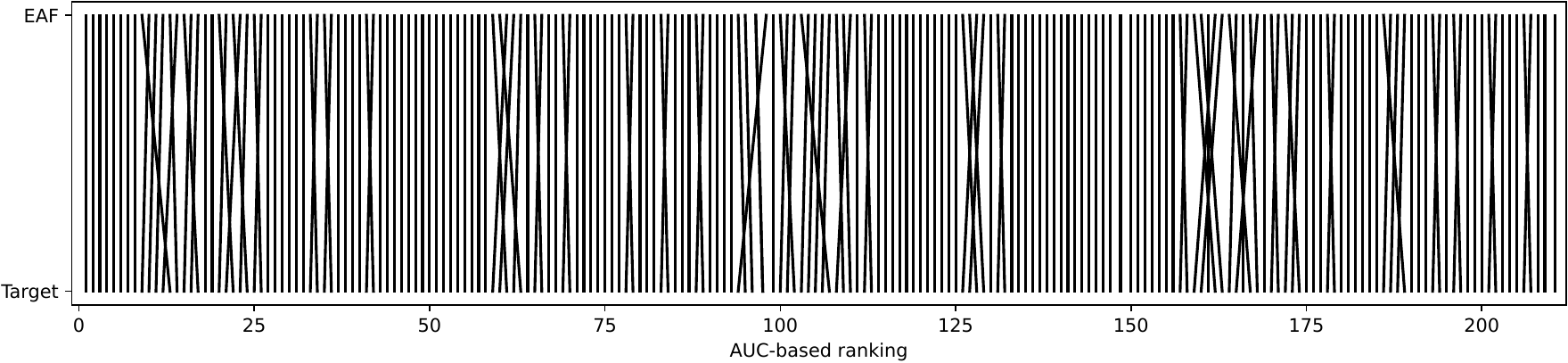}
    \caption{Ranking of 211 algorithms based on AUC (aggregated over the 24 10-dimensional BBOB problems) when using the EAF and the target-based ECDF with 51 targets. Lines crossing denote discrepancies in ranking. The sum of absolute rank differences (divided by two because each rank change
contributes twice when taking the absolute value) is 65. }
    \label{fig:rank_horizontal}
\end{figure*}

Finally, we may question whether these small differences matter at all when
ranking algorithms according to the AUC values. To answer this question, we
rank all 211 algorithms according to different AUC values, i.e., of the
EAF-based ECDF as well as the target-based ECDFs with 5, 10, 25 and 51
targets. We highlight the rank-differences between the 51 target and EAF setting in Figure~\ref{fig:rank_horizontal}.
We additionally calculate the rank differences between two rankings by summing the
absolute difference between ranks and dividing by two, because each rank change
contributes twice when taking the absolute value. In the case of the
2-dimensional BBOB problems, the rank differences between the target-based
ECDFs and the EAF-based ECDFs are: 298 (5 targets), 174 (10), 118 (25) and 101
(51). In the case of the 10-dimensional BBOB problems, the rank differences
between the target-based ECDFs and the EAF-based ECDFs are: 163 (5 targets),
81 (10), 49 (25), and 65 (51).
These numbers show that, although the differences in ranks quickly go down with
increasing number of targets, a high number of rank differences remain even with 51 targets.
Moreover, the last result shows that increasing the number of targets does not
always guarantee a more precise AUC value, since the precision depends on the
particular targets not only on their number.

\section{EAF-based Statistics}\label{sec:stats}

Several summary statistics can be computed starting from the EAF. In
particular, the level sets of the EAF are quantile-like functions that can be
represented as convergence curves. For example, the boundary region in the time
versus objective value space where the EAF has a value less than or equal to
0.5 is a synthetic convergence curve that represents the median convergence
curve (50\% percentile), such that any combination of time and objective
function value on or above this curve will be attained in 50\% of the runs of
an algorithm, but not necessarily all of them simultaneously in any single run. Similarly, the 100\%
percentile is the boundary that is attained by all runs of an algorithm and
represents the worst-case performance, whereas the lowest available percentile
represents the best-case performance. Examples of such percentiles are shown in
Figures~\ref{fig:eaf_bfgs} and~\ref{fig:eaf_cma}.

Similarly, considering the level sets of the EAF as quantile-like functions,
Binois et al.~\cite{BinGinRou2015gaupar} suggest the Vorob'ev expectation as a
useful definition of the expected value of a random
set~\cite[chap.~2]{Molchanov2005theory}. In a nutshell, the Vorob'ev
expectation is the quantile of the EAF with the closest AOCC value
to the expected value of the AOCC of a single run of the
algorithm. In other words, the Vorob'ev expectation is the synthetic
convergence curve that corresponds to a particular quantile of the EAF whose
AOCC value matches the mean AOCC of the actual runs of an algorithm. From a
benchmarking perspective, the Vorob'ev expectation is the ``mean'' convergence
curve and it can be used to summarize a large number of runs for
visualization and comparison purposes. In addition, dispersion statistics may be computed as
the probability of a single run deviating from this mean convergence
curve~\cite{BinGinRou2015gaupar}.

In addition to the above, the second-order EAF~\cite{FonGruPaq2005:emo}
measures the probability of attaining an objective function value $z$ not later
than $t$ if the algorithm has already obtained an objective value $z'$ not
later than $t'$. From a benchmarking perspective, such information is useful to
diagnose convergence issues, particularly, for algorithms that switch behavior
after a number of steps or evaluations.

Statistical measures and visualization techniques~\cite{ArzCebIru2022jcgs} that can be applied to the target-based ECDF are also applicable to the EAF-based ECDF. Some techniques may be directly applicable to the EAF. For example, it is possible to compute bootstrapped confidence bands~\cite{ArzCebIru2022jcgs} around the EAF that estimate the uncertainty. However, applying other techniques directly to the EAF may require additional research effort to handle the challenge of samples being sets of (runtime, quality) points rather than scalar numerical values.

\section{Conclusions}
\label{sec:conclusions}

 In this paper, we have argued that the EAF, and statistics derived from it, have significant benefits for the analysis and benchmarking of single-objective optimizers. In particular, we have shown how the widely-used target-based ECDF is just an approximation to the EAF-based ECDF and that the former eventually converges to the latter with increasing number of well-spread target values. We obtain a similar conclusion when considering the AUC of the target-based ECDF, which has been proposed as a measure of anytime performance. In this case, not only we argue that the AUC of the target-based ECDF is just an approximation of the AUC of the EAF-based ECDF, which is actually equal to the AUC of the EAF itself; but also we show that the latter is equivalent to the mean area over (or under in the case of maximization) the convergence  curves (AOCC) of the individual runs. Therefore, we conclude that  recording the mean AOCC value is sufficient for measuring anytime performance in situations where we are only interested in the AUC values, e.g., for benchmarking purposes or when guiding the tuning of parameters~\cite{LopStu2013ejor}.

In addition to the theoretical benefits of replacing the target-based ECDF with the EAF, we also demonstrate its practical use by means of integrating the eaf package into the IOHanalyzer platform. In particular, we argue that the visualization of EAF differences provides clear information about differences in anytime performance between optimizers.  We also show that the theoretical differences between the target-based ECDF and the EAF-based ECDF are measurable in real benchmarking data. Moreover, although the differences in value are small for the choice of targets recommended by the COCO platform, they still produce rank differences when comparing optimizers. Finally, the application of the EAF for analyzing single-objective optimizers opens the door to various EAF-based statistics that have not been explored in the benchmarking literature.

We have focused here on single-objective black-box algorithms for continuous optimization problems because that is the context where the target-based ECDF has gained popularity. However, our conclusions apply to benchmarking anytime algorithms in general, including gradient-based optimizers and optimizers for combinatorial optimization problems. 
Our conclusions also naturally extend to benchmarking multi-objective algorithms according to their anytime performance~\cite{RadLopStu2013emo,BroTusTusWag2016biobj}, where quality targets are given as vectors of objective values. However, the visualization~\cite{TusFil2011vis4d} of the EAF in three dimensions (i.e., 2 objectives plus runtime) and its computation in more than three dimensions~\cite{FonGueLopPaq2011emo} are still challenging.

\section*{Acknowledgment}
This work was initiated when M.\@ López-Ibáñez held an invited Professorship at Sorbonne University, France.
Part of this work was carried out during Dagstuhl Seminar 
23251 ``Challenges in Benchmarking Optimization Heuristics''
and at the ``Benchmarked: Optimization Meets Machine Learning'' workshop at the Lorentz center.

{\scriptsize
\bibliographystyle{IEEEtranN}
\providecommand{\MaxMinAntSystem}{{$\cal MAX$--$\cal MIN$} {Ant} {System}}
  \providecommand{\rpackage}[1]{{#1}}
  \providecommand{\softwarepackage}[1]{{#1}}
  \providecommand{\proglang}[1]{{#1}}

}

\vfill
\enlargethispage*{1.15cm}  
\begin{IEEEbiography}[{\includegraphics[width=1in,height=1.35in,clip,keepaspectratio]{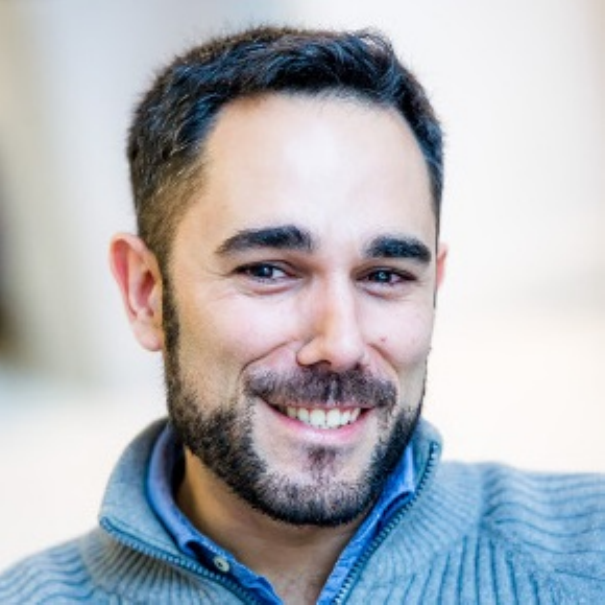}}]{Manuel L\'opez-Ib\'a\~nez}
   received his Master in Computer Science from the University of Granada, Spain in 2004 and his Ph.D. degree from Edinburgh Napier University, Edinburgh, U.K., in 2009. He is Full Professor at the Alliance Manchester Business School, University of Manchester, U.K.  He is  an elected member of the ACM SIGEVO Executive Board and the Business Committee of the Genetic and Evolutionary Computation Conference (GECCO), Editor-in-Chief of ACM Transactions on Evolutionary Learning and Optimization, Associate Editor of the Evolutionary Computation journal and Editorial Board Member of the Artificial Intelligence Journal.\  Prof.~López-Ibáñez has published more than 100 papers in international peer-reviewed journals and conferences on topics that include stochastic local search, black-box optimization, empirical reproducibility, multi-objective and interactive optimization algorithms for continuous and combinatorial problems, and the automatic configuration and design of optimization algorithms.
   \end{IEEEbiography}

\begin{IEEEbiography}[{\includegraphics[width=1in,height=1.35in,clip,keepaspectratio]{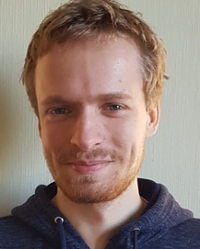}}]{Diederick Vermetten} is a PhD candidate at the Natural Computing group at LIACS, Leiden University. He is part of the core development team of the IOHprofiler project. His research interests include benchmarking of optimization heuristics, dynamic algorithm selection and configuration as well as hyperparameter optimization.
   \end{IEEEbiography}

\begin{IEEEbiography}[{\includegraphics[width=1in,height=1.25in,clip,keepaspectratio]{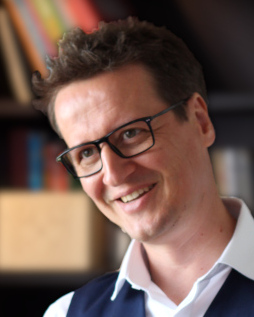}}]{Johann Dreo} is a senior research engineer at Institut Pasteur, working on decision support systems for helping bioinformatics against cancer. He received masters degrees in ecology in 1999, from Rennes~1 University, in applied mathematics in 2000, from Sorbonne University, and a Ph.D. in biomedical engineering from Paris-Sud University in 2003. He worked for several years as a researcher in Thales' industrial research center, on automated optimization solvers design. He is a member of the IEEE task force on Automated Algorithm Design, Configuration \& Selection, and the maintainer of the ParadisEO framework since 2010. His main interests are in hyperparameter optimization and explainable machine learning for computational biomedicine.
\end{IEEEbiography}

\begin{IEEEbiography}
[{\includegraphics[width=1in,height=1.25in,clip,keepaspectratio]{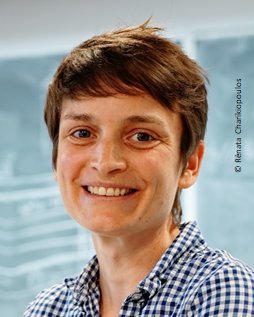}}]
{Carola Doerr} (formerly Winzen) received the Diploma degree in mathematics from Kiel University, Kiel, Germany, in 2007, the Dr.-Ing. degree from the Max Planck Institute for Informatics and Saarland University, Saarbrücken, Germany, in 2011, and the Habilitation degree (HDR) from Sorbonne Université, Paris, France, in 2020.
She is a CNRS Research Director of the LIP6 Computer Science Department, Sorbonne Université. From 2007 until 2012, she was a Business Consultant with McKinsey \& Company (on educational leave from 2010 onward). Her main research activities are in the analysis of black-box optimization algorithms.
Dr. Doerr is an Associate Editor of IEEE Transactions on Evolutionary Computation, ACM Transactions on Evolutionary Learning and Optimization, and Evolutionary Computation. She is a founding and a coordinating Member of the Benchmarking Network.
\end{IEEEbiography}

\end{document}